\documentclass[journal]{IEEEtran}

\IEEEoverridecommandlockouts                              % This command is only needed if
\usepackage{lineno}
\usepackage[cmex10]{amsmath}
\usepackage{amscd,amsfonts,amsmath,amssymb,mathrsfs,pifont,stmaryrd,tipa}
\usepackage{array,cases,dsfont,graphicx,texdraw}
\usepackage{graphics} % for pdf, bitmapped graphics files
\usepackage{epsfig} % for postscript graphics files
\usepackage{stfloats}
\usepackage{subfig}
\usepackage{hyperref}
\usepackage{epstopdf}
\modulolinenumbers[5]

\newtheorem{Remark}{Remark}

\newtheorem{Problem}{Problem}
\newenvironment{Proof}{\noindent{\em Proof:\/}}{\hfill $\Box$\par}
\newtheorem{Theorem}{Theorem}
\newtheorem{Lemma}{Lemma}

\newtheorem{Assumption}{Assumption}

\newcommand{\TAC}{{\it IEEE Transactions on Automatic Control}}

\newcommand{\R}{{\mathbb R}}
\newcommand{\EQQ}{\begin{eqnarray*}}
\newcommand{\ENN}{\end{eqnarray*}}
\newcommand{\EQ}{\begin{eqnarray}}
\newcommand{\EN}{\end{eqnarray}}
\title{\LARGE \bf
The Cooperative Output Regulation Problem of Discrete-Time Linear Multi-Agent Systems by the Adaptive Distributed Observer
\thanks{*This paper is the updated version of \cite{Huang17} where the definition of the function $\rho (A)$ in page 1 is revised so that it also applies to nonsymmetric matrices.}}

\author{Jie Huang% <-this % stops a space
\thanks{*This work has been supported  by
Hong Kong Special Administration Region, Research Grants Council, grant
No. 14200515.
}}% <-this % stops a space
\begin{document}

\maketitle
%\thispagestyle{empty}
%\pagestyle{empty}

%%%%%%%%%%%%%%%%%%%%%%%%%%%%%%%%%%%%%%%%%%%%%%%%%%%%%%%%%%%%%%%%%%%%%%%%%%%%%%%%
\begin{abstract}

In this paper, we first present an adaptive distributed observer for a discrete-time leader  system. This adaptive distributed  observer will provide, to each follower, not only the estimation of the leader's signal, but also the estimation of the leader's system matrix.
Then, based on the estimation of the matrix $S$, we  devise a discrete adaptive algorithm to calculate the solution to the regulator equations associated with each follower, and
obtain an estimated feedforward control gain.
Finally, we solve the cooperative output regulation problem for discrete-time linear multi-agent systems by both state feedback and output feedback adaptive distributed control laws utilizing the adaptive distributed observer.

\end{abstract}

%%%%%%%%%%%%%%%%%%%%%%%%%%%%%%%%%%%%%%%%%%%%%%%%%%%%%%%%%%%%%%%%%%%%%%%%%%%%%%%%
\section{Introduction}

The cooperative output regulation problem for continuous-time linear multi-agent systems using distributed observer approach was first studied  in \cite{shtac12}.
This problem aims to design a distributed
control law for a multi-agent system to achieve asymptotic tracking
of a class of reference inputs and rejection of a
class of disturbances.  Both reference inputs and disturbances are generated by a leader system.
The problem is an extension of the classical output regulation problem \cite{dtac76, fwa76, F} from a single plant to a multi-agent system.
On the other hand, the problem can also be viewed as a generalization of the leader-following consensus problem \cite{fmtac04, Hu1, jadbabaie2003,ni2010} because its objectives include not only asymptotic tracking but also disturbance rejection.

The core of the approach in \cite{shtac12} is the design of a distributed observer for the leader system of the following form,
\begin{equation}\label{2.2}
   \dot{v}=Sv
\end{equation}
where $v\in \R^{q}$ is the state of the leader system
representing the reference input to be tracked and/or the
external disturbance to be rejected
and $S\in \R^{q\times q}$ is a constant matrix. The distributed observer is capable of producing the estimation of the leader's signal to each follower
so that a distributed control law can be synthesized to solve the problem. However, a drawback associated with the distributed observer approach is that
 each follower needs to know the information of the matrix $S$. To remove this assumption, recently,  an adaptive distributed observer for the leader system
(\ref{2.2}) was proposed in \cite{ch15}.

The adaptive distributed observer not only estimates the leader's signal but also estimates the leader's dynamics. Thus, it does not
 require  every  follower know the matrix $S$.  In this paper, we will first propose a discrete counterpart of the adaptive distributed observer in \cite{ch15} for a discrete leader system of the following form
\begin{equation}\label{2.2d}
   v (t+1)=Sv (t).
\end{equation}
Then,  we will further develop an adaptive scheme to solve the cooperative output regulation problem for discrete linear multi-agent systems utilizing the discrete adaptive distributed observer.
Technically, we offer three specific contributions. Firstly, we
establish a stability result for a class of time-varying discrete-time systems that lends itself to the existence conditions of the adaptive distributed observer
of the leader system (\ref{2.2d}).
Secondly, based on the estimation of the matrix $S$, we  devise a discrete adaptive algorithm to calculate the solution to the regulator equations associated with each follower of the
discrete linear multi-agent system, and
obtain an estimated feedforward control gain.
Finally, we design both state feedback and output feedback  adaptive distributed control laws utilizing the discrete adaptive distributed observer to solve the cooperative output regulation problem for the discrete-time linear multi-agent system.

{\bf Notation.} $\otimes$ denotes
the Kronecker product of matrices.
%$1_N$ denotes an $N$ dimensional
%column vector whose components are all $1$.
$||x||$ denotes the
Euclidean norm of a vector $x \in \R^n$. $\mathbb{Z^+}$ denotes the set of nonnegative integers. Let $x: \mathbb{Z^+} \rightarrow \R^n$. Then we often denote $x (t)$, $ t \in  \mathbb{Z^+}$ by a shorthand notation $x$ when no confusion will occur.
$\sigma(A)$ denotes the spectrum of $A$ and  $\rho (A) =\max_{\lambda \in \sigma (A)} \{\frac{|\lambda|^2}{Re (\lambda)} \}$.
$\mathbf{1}_N$ denotes an $N\times1$ column vector whose elements are all $1$.
For $X_i\in \R^{n_i \times p}$, $i=1,\dots,m$,
col$(X_1,\dots,X_m)=[X_1^T,\dots,X_m^T]^T$. For any matrix
$A\in \R^{m\times n}$,
\begin{equation}
    \hbox{vec}(A)=\left[
                    \begin{array}{c}
                      A_1 \\
                      \vdots \\
                      A_n \\
                    \end{array}
                  \right]
\end{equation}
where $A_i$ is the $i$th column of $A$. For any column vector $X \in \R^{nq}$ for some positive integers $n$ and $q$,
\begin{equation}
    {M^q_n}(X)=\left[
                    \begin{array}{ccc}
                      X_1 \cdots  X_q \\
                    \end{array}
                  \right]
\end{equation}
where, for $i = 1, \dots, q$,  $X_i \in \R^n$, and  are such that $X = \mbox{col} (X_1, \dots, X_q)$.

\section{Problem Formulation}

Consider the following discrete-time linear multi-agent system:
\begin{subequations}\label{2.1}
    \begin{align}
      {x}_i (t+1) &=A_ix_i+B_i u_i+E_iv, ~ t  \in \mathbb{Z^+} \label{2.1.1}\\
      e_i (t) &=C_ix_i+D_i u_i+F_iv \label{2.1.2}\\
      y_{mi} (t)&=C_{mi}x_i+D_{mi}u_i+F_{mi}v, ~ i = 1, \dots, N \label{2.1.3}
      \end{align}
\end{subequations}
where $x_i\in \R^{n_i}$, $u_i\in \R^{m_i}$, $e_i\in \R^{m_i}$, $y_{mi}\in \R^{p_{mi}}$
are the state, control input, regulated output and measurement output of the
$i$th subsystem, respectively.

Like in \cite{shtac12}, we treat the system composed of (\ref{2.2d}) and  (\ref{2.1})
as a multi-agent system of $(N+1)$ agents with
(\ref{2.2d}) as the leader and the $N$ subsystems of (\ref{2.1}) as
$N$ followers, and define a graph\footnote{See \cite{let} or \cite{gr01} for a summary of digraph.}
$\bar{\mathcal{G}}=(\bar{\mathcal{V}},\bar{\mathcal{E}})$ with
$\bar{\mathcal{V}}=\{0,1,\dots,N\}$ and $\bar{\mathcal{E}}\subseteq
\bar{\mathcal{V}}\times \bar{\mathcal{V}}$. Here the node $0$ is
associated with the leader system (\ref{2.2d}) and the node $i$, $i =
1,\dots,N$, is associated with the $i$th subsystem of the
follower system (\ref{2.1}). For $i=0,1,\dots,N$, $j=1,\dots,N$,
$(i,j) \in \bar{\mathcal{E}}$ if and only if agent $j$ can use
the state or output of agent $i$ for control. Let $\bar{\mathcal{N}}_i$ denote
the neighbor set of the node $i$ of $\bar{\mathcal{G}}$.
% We can
%further define a subgraph $\mathcal{G}=(\mathcal{V},\mathcal{E})$ of
%$\bar{\mathcal{G}}$ where $\mathcal{V}=\{1,\dots,N\}$ and
%$\mathcal{E}\subseteq \mathcal{V}\times \mathcal{V}$ is obtained
%from $\bar{\mathcal{E}}$ by removing all the edges between node $0$
%and the nodes in $\mathcal{V}$.
%

We will consider the output feedback control law of the following form:
\begin{equation}\label{ctr22}
\begin{split}
u_i(t)=& k_i (z_i(t))\\
%\eta_i(t+1)=&g_2(\eta_i(t),\eta_j(t),j\in\bar{\mathcal{N}}_i)\\
z_i(t+1)=& g_i \big(z_i(t),   y_{mi}(t),  z_j(t),  y_{mj}(t), j\in\bar{\mathcal{N}}_i\big)\\
%\zeta_i(\theta)=& \zeta_{i0}(\theta), \theta\in I[-r,0]\\
i&=1,\cdots,N
\end{split}
\end{equation}
where $y_{m0} = v$,  $z_{i}\in \R^{l_i}$ for some integer $l_i$,  $k_i$ and $g_i$ are linear functions of their arguments whose specific form will be given in Section \ref{IV}.
It can be seen that, for each $i = 1, \cdots, N$, $j = 0, 1, \cdots, N$, $u_i$ of (\ref{ctr22}) depends on $y_{mj}$ only if the agent $j$ is a neighbor of the agent $i$. Thus,
the control law (\ref{ctr22}) is a distributed control law. The control law (\ref{ctr22}) contains the state feedback control law as a special case when $y_{mi} = x_i$, $i = 1, \cdots, N$.

We now define
the adaptive cooperative output regulation problem for \eqref{2.2d} and \eqref{2.1} as follows.
\begin{Problem}\label{prob1}
  Given systems \eqref{2.2d}, \eqref{2.1} and a graph $\bar{\mathcal{G}}$, design a distributed
  control law of the form (\ref{ctr22})
   %in the following form
%  \begin{subequations}\label{ctrlg}
%    \begin{align}
%      u_i&=f_i(\psi_i) \label{ctrlg.1}\\
%      \dot{\psi}_i&=k_i(\psi_i,\psi_j, j\in \bar{\mathcal{N}}_i), ~ i = 1, \dots, N \label{ctrlg.2}
%      \end{align}
%\end{subequations}
%  where $\psi_0$ represents the information of the leader system and $f_i$, $k_i$
%  are smooth functions,
  such that
  the trajectory of the closed-loop system starting from any initial state
  exists for all $t\geq 0$ and satisfies
  \begin{itemize}
    \item when $v$ is bounded, the trajectory of the closed-loop system is bounded for all $t\geq 0$;
    \item the regulated output satisfies
    \begin{equation}
    \lim_{t\rightarrow\infty}e_i(t)=0,\ i=1,\dots,N.
  \end{equation}
  \end{itemize}
\end{Problem}

We need the following assumptions for the solvability of Problem \ref{prob1}.

\begin{Assumption}\label{ass1}
 All the eigenvalues of $S$ have modulus smaller than or equal to $1$.
\end{Assumption}

\begin{Assumption}\label{ass2}
  For $i=1,\dots,N$, $(A_i,B_i)$ are stabilizable.
\end{Assumption}

\begin{Assumption}\label{ass2.1}
  For $i=1,\dots,N$, $(C_{mi},A_i)$ are detectable.
\end{Assumption}

\begin{Assumption}\label{ass3}
For $i=1,\dots,N$,   the following linear matrix equations
  \begin{subequations}\label{re}
    \begin{align}
      X_i S &= A_i X_i+B_i U_i+E_i \label{re.1}\\
      0 &= C_i X_i+D_i U_i+F_i \label{re.2}
      \end{align}
      \end{subequations}
      have  unique solution pairs $(X_i, U_i)$.
\end{Assumption}
\begin{Assumption}\label{ass4}
  The graph $\bar{\mathcal{G}}$ contains a spanning
  tree with the node $0$ as the root.
\end{Assumption}
\begin{Remark}

A large class of signals such as the step function, ramp function, and sinusoidal function satisfies Assumption \ref{ass1}.

%%Moreover, if the  cooperative output regulation problem  is solvable under Assumption \ref{ass1}, then the same problem is also solvable even if $S$ contains some eigenvalues with modulus smaller than unity because the modes corresponding to the eigenvalues of $S$ with modulus smaller than 1 will decay to zero asymptotically.

In the classical linear output regulation problem, equations \eqref{re} are
  called the regulator equations whose solvability imposes a necessary
  condition for the solvability of the output regulation problem.
By Theorem 1.9 of \cite{h04}, for any matrices $E_i$ and $F_i$,  the regulator equations
(\ref{re}) are solvable if and only if  \EQ \mbox{\rm    rank}
\left [
\begin{array}{cc}
A_i - \lambda I_{n_i} & B_i \\
C_i & D_i
\end{array} \right ]
= n_i+m_i,  ~~\forall~~\lambda \in \sigma(S). \label{tzeros} \EN

%
%
% \cite{dtac76,fco77,fwa76}.
%  The solution pair $(X_i,U_i)$ determines the feedforward control gains
%  for both the state feedback and measurement output feedback controllers.
%
%
%
%  However, in our case, since not all the followers can access the
%  system matrix $S$, we cannot calculate $(X_i,U_i)$, and hence the feedforward gains in advance.
%  To overcome this difficulty, we will develop an auxiliary system to solve
%  the solutions of \eqref{re} online.
Assumption \ref{ass4} is a standard assumption in the literature of the cooperative control of multi-agent systems subject to static networks.
\end{Remark}

\section{Adaptive Distributed Observer}

 The distributed observer for a discrete leader system of the form (\ref{2.2d}) was proposed in \cite{let} and takes the following form:
\begin{equation}\label{do.2x}
    {\eta}_i (t+1) =S\eta_i (t) +\mu S \sum_{j=0}^Na_{ij}(\eta_j (t) -\eta_i (t)), i=1,\dots,N
\end{equation}
where $\eta_0=v$, $\eta_i\in \R^q$, $i = 1, \dots, N$, $\mu>0$, and $\bar{\mathcal{A}} =[a_{ij} ]_{i,j =0}^{N}$ denote the weighted adjacency matrix of $\bar{\mathcal{G}}$.
Let $\eta(t)= \mbox{col} (\eta_1(t),\cdots, \eta_N(t))$, $\hat{v}(t)=\mathbf{1}_N\otimes v(t)$, and $\tilde{\eta}(t)=\eta(t)-\hat{v} (t)$. Then (\ref{do.2x}) can be put in the following form:
\begin{equation}\label{obt}
\begin{split}
\tilde{\eta} (t+1)&=\big( (I_N\otimes S)-\mu(H\otimes S)\big)\tilde{\eta}(t)\\
\end{split}
\end{equation}
where $H=[h_{ij}]\in \R^{N\times N}$ with $h_{ii}=\sum_{j=0}^Na_{ij}$ and $h_{ij}=-a_{ij}$, for any $i\neq j$. By Lemma 1 of \cite{Hu1},
under Assumption \ref{ass4}, all the eigenvalues of $H$ have positive real parts.

If there exists some $\mu$ such that the matrix $\big( (I_N\otimes S)-\mu(H\otimes S) \big )$ is Schur, then, for any $v
(0)$, and $\eta_i (0)$, $ i = 1, \dots, N$, we have \EQ \label{p0}
\lim_{t \rightarrow \infty} (\eta_i(t) - v (t)) = 0. \EN
Thus, we call the system \eqref{do.2x} a distributed observer of the leader system  if and only if the system \eqref{obt} is asymptotically stable.
A detailed discussion on the stability of the system \eqref{obt}  is summarized in Lemma 3.1 of \cite{let}. In particular,
denote the eigenvalues of $S$ by $\{\lambda_1, \cdots, \lambda_q\}$ where $0\leq |\lambda_1| \leq \cdots \leq |\lambda_q|$, and the eigenvalues of $H$ by $\{a_l\pm jb_l\}$, where $b_l=0$ when $1\leq l \leq N_1$ with $0 \leq N_1 \leq N$ and $b_l\neq 0$ when $(N_1+1)\leq l\leq N_2$ where $N_1+2(N_2-N_1)=N$.  Then,  under Assumptions \ref{ass1} and \ref{ass4},
the matrix $((I_N\otimes S-\mu(H\otimes S))$ is Schur for all $\mu$ satisfying
\EQ \label{conmu}
 \max_{l = 1, \cdots, N} \bigg \{ \frac{ a_l -\sqrt{\Delta_l }}{a_l^2+b_l^2} \bigg \} < \mu  < \min_{l = 1, \cdots, N} \bigg \{ \frac{ a_l  + \sqrt{\Delta_l}}{a_l^2+b_l^2} \bigg \}
\EN
where $\Delta_l = \frac{(a_l^2+b_l^2)}{ |\lambda_l|^2 } -  b^2_{l}$. If $|\lambda_l| = 1$ for all $l = 1, \cdots, q$, then $\Delta_l =a_l^2$ and (\ref{conmu}) reduces to
\EQ \label{conmu2}
0  < \mu  < \min_{l = 1, \cdots, N} \bigg \{ \frac{ 2 a_l }{a_l^2+b_l^2} \bigg \}.
\EN
%\EQ \label{conmu}
% \max_{l = 1, \cdots, N} \bigg \{ \frac{|\lambda_q| -1 }{a_l |\lambda_q| } \bigg \} < \mu  < \min_{l = 1, \cdots, N} \bigg \{ \frac{|\lambda_q| +1 }{a_l |\lambda_q| }  \bigg \}.
%\EN
%Under Assumption \ref{ass4}  and the additional assumption that the digraph $\mathcal{\bar{G}}$  is undirected,
%the matrix $(I_N\otimes S-\mu(H\otimes S))$ is Schur if and only if
%  \EQ \label{mu1x}
%\max_{l = 1, \cdots, N} \bigg \{ \frac{|\lambda_q| -1}{a_l |\lambda_q|}  \bigg \} <   \min_{l = 1, \cdots, N}  \bigg \{  \frac{|\lambda_q| + 1}{a_l |\lambda_q|} \bigg \}.
%\EN
%And, if the condition (\ref{mu1x}) is satisfied, the matrix $(I_N\otimes S_0-\mu(H\otimes S_0))$ is Schur for all $\mu$ satisfying
%\EQ \label{mu1xx}
%\max_{l = 1, \cdots, N} \bigg \{ \frac{|\lambda_q| -1}{a_l |\lambda_q|}  \bigg \} <  \mu < \min_{l = 1, \cdots, N}  \bigg \{  \frac{|\lambda_q| + 1}{a_l |\lambda_q|} \bigg \}.
%\EN

However,  in (\ref{do.2x}),  the matrix $S$ is used by every follower, which may not be realistic in some  applications. Here, we will further propose the following so-called  adaptive distributed  observer candidate:
\begin{subequations}\label{do}
    \begin{align}
      {S}_i (t+1)&= {S}_i + \mu_1\sum_{j=0}^Na_{ij}(S_j-S_i) \label{do.1}\\
      {\eta}_i (t+1) &= S_i\eta_i+\mu_2 S_i \sum_{j=0}^Na_{ij}(\eta_j-\eta_i)  \label{do.2}
      \end{align}
\end{subequations}
where $S_0=S$, $S_i\in \R^{q\times q}$, $i=1,\dots,N$, $\mu_1,\mu_2>0$.
If there exist some $\mu_1,\mu_2>0$ such that, for $i=1,\dots,N$, for any initial condition, the solution to (\ref{do})
satisfies  $\lim_{t\rightarrow\infty}(S_i(t)- S)=0$, $\lim_{t\rightarrow\infty}({\eta}_i(t)- v (t))=0$, then
  (\ref{do}) is called the adaptive  distributed observer for the leader system.
It can be seen from (\ref{do}) that only those followers who are the children of the leader know the matrix $S$. Thus,
 the adaptive distributed  observer is more realistic than the distributed observer.

To find the conditions under which (\ref{do}) is an adaptive  distributed observer for the leader system, we
first establish the following lemma.

\begin{Lemma}\label{lem0}
  Consider the following system
  \begin{equation}\label{lem0.1}
    x (t+1)=F x (t)+F_1(t)x (t)+F_2(t)
  \end{equation}
  where $x\in \R^n$, $F\in \R^{n\times n}$ is Schur, $F_1(t)$ and $F_2(t)$
  are well defined  for all $t\in \mathbb{Z^+}$. If
  $F_1(t),F_2(t)\rightarrow0$  (exponentially) as $t\rightarrow\infty$,
  then, for any $x (0) \in\R^n$,  $x(t)\rightarrow0$ (exponentially) as $t\rightarrow\infty$.
\end{Lemma}

\begin{Proof}
  If $F_2(t)=0$, then \eqref{lem0.1} reduces to
  \begin{equation}\label{lem0.2}
    x (t+1) =F x (t)+F_1(t)x (t).
  \end{equation}
  Since $F$ is Schur, for any symmetric and positive definite matrix $N\in \R^{n\times n}$, there exist a symmetric positive definite matrix $M \in \R^{n\times n}$ such that $F^T M F - M = - N$. Let $V (t) = x^T (t) M x (t)$. Then, along the trajectory of \eqref{lem0.2}
    \EQQ \label{lem0.3}
    &&V (t+1) - V (t) \\
    &=& (F x (t)+F_1(t)x (t))^T M (F x (t)+F_1(t)x (t))\\
    & - &x^T (t) M x (t)\\
    &=& x^T (t) F^T M F x (t) + x^T (t) (2 F^T M F_1(t)\\
    & +& F^T_1(t) M F_1(t)) x (t) - x^T (t) M x (t)  \\
      &=& - x^T (t) (N - 2 F^T M F_1(t) - F^T_1(t) M F_1(t)) x (t)
  \ENN
  Since $F_1(t)\rightarrow0$ as $t\rightarrow\infty$,  there exist $T  \in \mathbb{Z^+}$ and $\lambda > 0$ such that $x^T (t) (N - 2 F^T M F_1(t) - F^T_1(t) M F_1(t)) x (t) > \lambda
  || x (t)||^2$ for all $t \geq T$. Thus, the system (\ref{lem0.2}) is exponentially stable.
  Therefore, for any initial condition, the solution to \eqref{lem0.1} is bounded for all $t\in \mathbb{Z^+}$.
  Since $F$ is Schur,  system \eqref{lem0.1} is input-to-state stable
  with $F_1(t)x + F_2(t)$ as input. Thus, by Lemma 3.8 of \cite{jiang01},  the system \eqref{lem0.1}   has  ${\cal K}$ asymptotic gain property, i.e.,  there exists   a class $\mathcal{K}$ function
  $\gamma$ such that, for any $x(0)\in \R^n$, the solution $x(t)$ of \eqref{lem0.1} satisfies
  \begin{equation*}
  { \limsup_{t \rightarrow\infty}}||x(t)|| \leq \gamma (\limsup_{t\rightarrow\infty} ||F_1 (t) x(t) + F_2(t)||).
  \end{equation*}
 Since $x(t)$ is bounded, if  $F_1(t),F_2(t)$ tend to zero (exponentially), so does $x (t)$.
 \end{Proof}

%Let $\bar{\mathcal{A}}=[a_{ij}]$ be the adjacency matrix of $\bar{\mathcal{G}}$.
%Since $S$ is not accessible to all the followers, for each follower, we
%design the following distributed compensator
%\begin{equation}\label{de}
%    \dot{S}_i=\mu_1\left(\sum_{j=1}^Na_{ij}(S_j-S_i)+a_{i0}(S-S_i)\right)
%\end{equation}
%where $S_i\in R^{q\times q}$, $\mu_1>0$. \eqref{de} is called the distributed
%estimator for the leader system \eqref{2.2} in the sense of the following lemma

We now establish the main result of this section.

\begin{Lemma} \label{lemma1}
  Given the systems \eqref{2.2d} and \eqref{do}, let $\tilde{S}_i=S_i-S$,
  $\tilde{\eta}_i = {\eta}_i  - v  $.
 Then, for any  $S_i(0)$ and $\eta_i(0)$,
  we have \\
 (i) Under Assumption \ref{ass4}, for any $\mu_1$ satisfying $ 0 < \mu_1 < \frac{2}{\rho (H)}$, for $i=1,\dots,N$,
 \begin{equation}\label{nn.1}
    \lim_{t\rightarrow\infty}\tilde{S}_i(t)=0
  \end{equation}
  exponentially, and \\
(ii)  Under Assumptions \ref{ass1} and \ref{ass4}, let $\mu_1$ satisfy $0 < \mu_1 < \frac{2}{\rho (H)}$ and
  let $\mu_2$ be such that the matrix $(I_N\otimes S)-\mu_2(H\otimes S)$ is
Schur. Then, for $i=1,\dots,N$, for any $\tilde{\eta}_i(0)$,
  \begin{equation}\label{sc}
    \lim_{t\rightarrow\infty}\tilde{\eta}_i(t)=0
  \end{equation}
  exponentially.
\end{Lemma}
\begin{Proof}
Part (i). Let $ \tilde{S} = \mbox{col} ({\tilde{S}}_1, \dots, {\tilde{S}}_N)$.
Then \eqref{do.1} can be put in the following form
\EQ
{\tilde{S}} (t+1) = (I_{Nq}- \mu_1 (H \otimes I_q)) \tilde{S} (t)
\EN
%or a more standard form
%\begin{equation}
%    \hbox{vec}(\tilde{S}  (t+1)) = ( (I_{Nq^2} - \mu_1 (I_q \otimes H \otimes I_q))\hbox{vec}(\tilde{S} (t)).
%\end{equation}
Under Assumption \ref{ass4}, by Lemma 1 of \cite{Hu1},
  all the eigenvalues of $H$ have positive real parts. Thus, for any $\mu_1$ satisfying $ 0 < \mu_1 < \frac{2}{\rho (H)}$, the matrix $  ( (I_{Nq} - \mu_1 ( H \otimes I_q))$ is Schur.
 Therefore, $\lim_{t\rightarrow\infty} \tilde{S} (t)=0$ exponentially, that is,  for $i=1,\dots,N$, $\lim_{t\rightarrow\infty}\tilde{S}_i(t)=0$ exponentially.

Part (ii).  By \eqref{do}, we have
  \begin{equation}\label{lem3.1}
  \begin{aligned}
   {\tilde{\eta}}_i (t+1) &=S_i\eta_i-Sv+\mu_2  S_i \sum_{j=0}^Na_{ij}
    (\tilde{\eta}_j-\tilde{\eta}_i)\\
    %&=S_i\eta_i-S\eta_i+S\eta_i-Sv+\mu_2 S_i \sum_{j=0}^Na_{ij}
   % (\tilde{\eta}_j-\tilde{\eta}_i)\\
    &=S\tilde{\eta}_i+\tilde{S}_i\eta_i+\mu_2 S_i \sum_{j=0}^Na_{ij}
    (\tilde{\eta}_j-\tilde{\eta}_i) \\
    &=S\tilde{\eta}_i+\tilde{S}_iv+\tilde{S}_i\tilde{\eta}_i+\mu_2 S_i \sum_{j=0}^Na_{ij}
    (\tilde{\eta}_j-\tilde{\eta}_i) \\
    &=S\tilde{\eta}_i+ \mu_2 S \sum_{j=0}^Na_{ij}
    (\tilde{\eta}_j-\tilde{\eta}_i) \\
    & + \tilde{S}_iv + \tilde{S}_i\tilde{\eta}_i +\mu_2 \tilde{S}_i \sum_{j=0}^Na_{ij}
    (\tilde{\eta}_j-\tilde{\eta}_i).
  \end{aligned}
  \end{equation}
 Let $ \tilde{\eta}= \mbox{col} (\tilde{\eta}_1, \dots, \tilde{\eta}_N)$ and $\tilde{S}_d =\hbox{block diag}\{\tilde{S}_1,\dots,\tilde{S}_N\}$. Then, \eqref{lem3.1}
 can be rewritten in the following compact form
\begin{equation}
 \begin{aligned}
& {\tilde{\eta}} (t+1)  =((I_N\otimes S)-\mu_2(H\otimes S))\tilde{\eta}
    +\tilde{S}_d ( \mathbf{1}_N \otimes v)  \\
& + \left ( \tilde{S}_d - \mu_2\left [ \begin{array}{c} H_1 \otimes \tilde{S}_1 \\
\vdots \\
H_N \otimes \tilde{S}_N \end{array}
\right ] \right )\tilde{\eta}
\end{aligned}
\end{equation}
where, for $i = 1, \cdots, N$, $H_i$ is the $i{th}$ row of $H$. By assumption, the matrix $((I_N\otimes S)-\mu_2(H\otimes S))$ is
Schur. By Part (i) of this Lemma,  $\lim_{t\rightarrow\infty}\tilde{S}(t)=0$
exponentially, Thus,  under Assumption \ref{ass1},
$\tilde{S}_d (t) (\mathbf{1}_N \otimes v(t))$ will decay to
zero exponentially, too. It follows from Lemma \ref{lem0} that, for any $\tilde{\eta}_i(0)$,
\begin{equation}
     \lim_{t\rightarrow\infty}\tilde{\eta}_i(t)=0
\end{equation}
exponentially and the proof is completed.
\end{Proof}
%
%\begin{Remark}
%v unbounded
%\end{Remark}

\section{Main Result}\label{IV}

When the matrix $S$ is known by every follower, a control law utilizing the  solution to the regulator equations has been designed in \cite{let} for solving our problem.
Since, in this paper, we assume  those followers which are not the children of the leader do not know $S$, we cannot directly use the  solution to the regulator equations
to design our control law. We will propose to adaptively calculate the solution to the regulator equations based on the
estimation $S_i$ of $S$.
For this purpose,  we need to establish the following lemma.

%\subsection{A Lemma}
%
% First, we introduce the following lemma.
\begin{Lemma}\label{lem2}
 Consider the following linear algebraic equation:
  \begin{equation}\label{lae1}
   A x = b
  \end{equation}
 where $x\in \R^n$, $A\in \R^{n\times n}$ is nonsingular, and  $b\in \R^n$.
 Let $B (t)\in \R^{n\times n}$ be well defined for all $t\in \mathbb{Z^+}$ such that
  $\tilde{A}(t)\triangleq B(t)-A\rightarrow0$ exponentially
  as $t\rightarrow\infty$.
  Then,  for any $x(0)\in \R^n$ and $ 0 < \varepsilon < \frac{2}{\rho (A^T A)}$, the solution $x (t)$  to the following system
  \begin{equation}\label{lem2.2}
   x (t+1)= x (t)-\varepsilon B(t)^T(B(t)x (t)-b)
  \end{equation}
 is such that
  \begin{equation}
    \lim_{t\rightarrow\infty}(x(t)-A^{-1} b)=0
  \end{equation}
  exponentially.
\end{Lemma}

\begin{Proof}
\begin{equation}\label{lem2.5}
    \begin{aligned}
    {x} (t+1)  &= x  -\varepsilon ( B(t)^T B(t)x - B(t)^Tb)\\
      &= {x}  -\varepsilon B(t)^T B(t)x+\varepsilon B(t)^Tb\\
      &= {x}  -\varepsilon A^TAx+\varepsilon A^TAx-\varepsilon B(t)^T B(t)x\\
      &+\varepsilon B(t)^Tb-\varepsilon A^Tb+\varepsilon A^Tb\\
      &={x} -\varepsilon A^TA {x}+\varepsilon (A^TA-B(t)^TB(t))x\\
      &+\varepsilon \tilde{A}(t)^Tb+\varepsilon A^Tb\\
      &= {x}  -\varepsilon A^TA {x}+\varepsilon A^Tb  +  F (t) {x}  +d(t).
    \end{aligned}
  \end{equation}
  where
  \begin{equation}
   F (t) = \varepsilon (A^TA-B(t)^TB(t)),   d(t)= \varepsilon \tilde{A}(t)^Tb.
  \end{equation}

Let $\bar{x}={x}-{x}^*$ where $ x^*= A^{-1} b$. Then, from (\ref{lem2.5}),  we have
\begin{equation}
\begin{aligned}
    {\bar{x}} (t+1)&= \bar{x} (t) -\varepsilon A^TA (\bar{x} (t) + {x}^*) + \varepsilon A^Tb \\
      &+  F (t) (\bar{x} (t) + {x}^*)  +d(t) \\
      &= (I_n -\varepsilon A^TA ) \bar{x} (t)  \\
      &+  F (t) \bar{x} (t) +  (F (t) {x}^*+d(t)).
    \end{aligned}
\end{equation}

Since  $\lim_{t\rightarrow\infty}\tilde{A}(t)=0$
exponentially, both $F (t)$ and
$d (t)$ will decay to zero
exponentially. Also, since our choice of $\varepsilon$ is such that $(I_n -\varepsilon A^TA )$ is Schur,  by Lemma \ref{lem0},
for any $x(0)\in \R^n$, \begin{equation}
 \lim_{t\rightarrow\infty}(x(t)-x^*) = 0
\end{equation}
exponentially. Therefore, the proof is completed.
\end{Proof}

%By letting $A(t)=A$ for all $t\geq 0$ in Lemma \ref{lem2}, we have the following result.
%\begin{Corollary}\label{cor}
%Let $A\in R^{m\times n}$,  $b\in R^m$, and  $\hbox{rank}(A)= \hbox{rank}(A, b)= k $ for
%some positive integers $m,n,k\geq 1$.
%Then, for any $x(0)\in R^n$ and $\epsilon>0$,   the following system
%  \begin{equation}\label{lem1.2}
%    \dot{x}=-\epsilon A^T(Ax-b)
%  \end{equation}
%has a unique bounded solution $x (t)$ for $t \geq 0$ such that,  for some $x^* \in R^n$ satisfying $
%    Ax^*=b$,
% \begin{equation}
%    \lim_{t\rightarrow\infty}(x(t)-x^*)=0
%  \end{equation}
%  exponentially.
%\end{Corollary}

\begin{Remark}
 Using a differential equation to solve a linear algebraic equation of the form (\ref{lae1})  was studied in \cite{cu92} for the case where $A$ is known, and
 was  studied in  \cite{ch15} for the case
 where $A$ is unknown. Lemma \ref{lem2} here further shows how to solve (\ref{lae1}) with $A$ unknown
by the difference equation (\ref{lem2.2}).
\end{Remark}

%\subsection{Adaptive solution of the regulator equations}

Using the information of $S_i$, an adaptive algorithm was proposed in \cite{ch15} to calculate the solution to the regulator equations.
This algorithm is governed by a set of nonlinear differential equations. Here we will develop a discrete counterpart of the adaptive algorithm in \cite{ch15}
to calculate the solution to the regulator equations by a set of difference equations.
For this purpose, like in \cite{ch15}, for $i=1,\dots,N$, let
\begin{equation}
    x_i=\hbox{vec}\left(\left[
             \begin{array}{c}
               X_i \\
               U_i \\
             \end{array}
           \right]\right),\ b_i=\hbox{vec}\left(\left[
             \begin{array}{c}
               E_i \\
               F_i \\
             \end{array}
           \right]\right),
\end{equation}
\begin{equation}
    Q_i=S^T\otimes \left[
      \begin{array}{cc}
        I_{n_i} & 0 \\
        0 & 0 \\
      \end{array}
    \right]-I_q\otimes \left[
                      \begin{array}{cc}
                        A_i & B_i \\
                        C_i & D_i \\
                      \end{array}
                    \right],
\end{equation}
and
\begin{equation}
    G_i(t)=S_i(t)^T\otimes \left[
      \begin{array}{cc}
        I_{n_i} & 0 \\
        0 & 0 \\
      \end{array}
    \right]-I_q\otimes \left[
                      \begin{array}{cc}
                        A_i & B_i \\
                        C_i & D_i \\
                      \end{array}
                    \right]
\end{equation}
where $S_i (t)$ is generated by (\ref{do.1}). It is noted that the dimensions of $x_i$, $Q_i$ and $G_i (t)$ are $q(n_i+m_i)  \times 1$,  $q (n_i+m_i) \times q (n_i+m_i)$ and $q (n_i+m_i) \times q (n_i+m_i)$, respectively.

We have the following result.

\begin{Lemma}\label{lem4}
Under Assumption \ref{ass3},  for $i=1,\dots,N$, for any initial condition $\zeta_i (0)$,
each of the following equations
\begin{equation}\label{ctrlols}
    {\zeta}_i (t+1) = {\zeta}_i (t)-\mu_{3i} G_i(t)^T(G_i(t)\zeta_i (t)-b_i)
\end{equation}
where $0 < \mu_{3i} < \frac{2}{\rho (Q_i^T Q_i)}$, $i = 1, \cdots, N$, has a  bounded solution for all $t\in \mathbb{Z^+}$.
Moreover, let $\Xi_i(t) = M_{(n_i+m_i)}^q (\zeta_i (t))$.
Then,
\begin{equation} \label{matrixs}
    \lim_{t\rightarrow\infty} \left ({\Xi}_i(t) - \left[ \begin{array}{c}
               X_i\\
               U_i \\
             \end{array}  \right] \right) =0
\end{equation}
exponentially.
\end{Lemma}

\begin{Proof}
The regulator equations \eqref{re} can be put in the
following form
\begin{equation}\label{rec}
    \left[
      \begin{array}{cc}
        I_{n_i} & 0 \\
        0 & 0 \\
      \end{array}
    \right]\left[
             \begin{array}{c}
               X_i \\
               U_i \\
             \end{array}
           \right]S-\left[
                      \begin{array}{cc}
                        A_i & B_i \\
                        C_i & D_i \\
                      \end{array}
                    \right]\left[
             \begin{array}{c}
               X_i \\
               U_i \\
             \end{array}
           \right]=\left[
                     \begin{array}{c}
                       E_i \\
                       F_i \\
                     \end{array}
                   \right].
\end{equation}
By Theorem 1.9 of \cite{h04},   \eqref{rec}
can be transformed into the following form
\begin{equation}\label{cre}
    Q_ix_i=b_i.
\end{equation}
Moreover, Assumption \ref{ass3} holds if and only if $Q_i$ is nonsingular for $i = 1, \dots, N$. Thus,
for $i = 1, \dots, N$, the linear algebraic equations (\ref{cre}) has a unique solution $Q_i^{-1} b_i$.
By \eqref{nn.1}, $\lim_{t\rightarrow\infty}(G_i(t)-Q_i)=0$ exponentially.
Therefore, by Lemma \ref{lem2}, for any ${\zeta}_i (0)$,
the solution $\zeta_i(t)$  to (\ref{ctrlols}) is  such that
\begin{equation} \label{vectors}
    \lim_{t\rightarrow\infty}(\zeta_i(t)-Q_i^{-1} b_i)=0
\end{equation}
exponentially.  (\ref{vectors}) implies (\ref{matrixs})  since
\begin{equation}
\left (\Xi_i (t) -   \left[  \begin{array}{c}
             X_i \\
               U_i
             \end{array} \right] \right ) = M_{(n_i+m_i)}^q (\zeta_i(t)-Q_i^{-1} b_i).
\end{equation}
\end{Proof}

\begin{Remark}
It is noted that, for $i = 1, \cdots, N$,  $\mu_{3i}$ in (\ref{ctrlols}) depends on the matrix $S$. Since
$\mu_{3i}$ can be calculated off-line, the algorithm (\ref{ctrlols}) itself does not need to know $S$ once $\mu_{3i}$ has been predetermined.
\end{Remark}

%\subsection{Solvability of the Problem}
We now ready to show how to solve Problem \ref{prob1}
by state feedback control. Partition $\Xi_i (t)$ as $\Xi_i (t)=[X_i(t)^T,U_i(t)^T]^T$, where $X_i(t)\in \R^{n_i\times q}$ and $U_i(t)\in \R^{m_i\times q}$.
Since $(A_i,B_i)$ is stabilizable, let $K_{xi}$ be such that
$\tilde{A}_i\triangleq A_i+B_iK_{xi}$ is Schur and $K_{\eta i}(t)$ be given as
\begin{equation} \label{keta}
    K_{\eta i}(t)=U_i(t)-K_{xi}X_i(t).
\end{equation}
For $i=1,\dots,N$, we design the following state feedback controller
\begin{equation}\label{ctrlu}
    u_i=K_{xi}x_i+K_{\eta i}(t)\eta_i.
\end{equation}

We have the following result.
\begin{Theorem}\label{thm1}
Under Assumptions \ref{ass1}, \ref{ass2}, \ref{ass3} and \ref{ass4},  let $\mu_1$ satisfy $0 < \mu_1 < \frac{2}{\rho (H)}$, $\mu_2$ be such that the matrix $(I_N\otimes S)-\mu_2(H\otimes S)$ is
Schur, and $0 < \mu_{3i} < \frac{2}{\rho (Q_i^T Q_i)}$, $i = 1, \cdots, N$. Then, Problem \ref{prob1} is solvable by the control law composed of  \eqref{do}, \eqref{ctrlols}
  and \eqref{ctrlu}.
\end{Theorem}
\begin{Proof}
Let $\tilde{x}_i (t) =x_i (t) -X_{i}v (t)$, $\tilde{u}_i (t) =u_i (t) -U_{i} v (t)$,
$K_{\eta i}=U_i-K_{xi}X_i$, and $\tilde{K}_{\eta i} (t) = {K}_{\eta i} (t) - K_{\eta i}$.
By making use of the solution to the regulator equations \eqref{re}, we obtain, for $i=1,\dots,N$,
\begin{equation}\label{sob0.1}
   \begin{aligned}
    {\tilde{x}}_i (t+1) &=A_ix_i+B_iu_i+E_iv-X_iSv\\
     &=A_i(\tilde{x}_i+X_iv)+B_i(\tilde{u}_i+U_iv)+E_iv-X_iSv\\
     &=A_i \tilde{x}_i + B_i \tilde{u}_i
   \end{aligned}
\end{equation}
and
\begin{equation}\label{sob0.2}
   \begin{aligned}
     e_i (t) &=C_ix_i+D_iu_i+F_iv\\
     &=C_i(\tilde{x}_i+X_i v)+D_i(\tilde{u}_i+U_i v)+F_iv\\
     &=C_i \tilde{x}_i + D_i \tilde{u}_i.
   \end{aligned}
\end{equation}
Further, we have
\begin{equation}\label{sob2}
    \begin{aligned}
      \tilde{u}_i (t) &= K_{xi} (\tilde{x}_i + X_{i} v) + K_{\eta i} (t) ( \tilde{\eta}_i  + v) - U_{i} v \\
   &=K_{xi}\tilde{x}_i + (U_i-K_{\eta i})v+ K_{\eta i} (t) ( \tilde{\eta}_i  + v) - U_{i} v\\
   &=K_{xi}\tilde{x}_i + K_{\eta i}(t)\tilde{\eta}_i+\tilde{K}_{\eta i}(t)v.
    \end{aligned}
\end{equation}
Substituting \eqref{sob2} to \eqref{sob0.1} gives
\begin{equation}\label{sob3}
   {\tilde{x}}_i (t+1) = (A_i+B_iK_{xi})\tilde{x}_i + f_i (t)
  \end{equation}
where $f_i(t)=B_iK_{\eta i}(t)\tilde{\eta}_i+B_i\tilde{K}_{\eta i}(t)v$.
By Lemma \ref{lemma1}, $\tilde{\eta}_i (t)$
decays to zero exponentially.
Since $\tilde{K}_{\eta i} (t) = (U_i(t)-U_i)-K_{xi}(X_i(t)-X_i)$, by Lemma \ref{lem4},
$\lim_{t \rightarrow \infty} \tilde{K}_{\eta i} (t) =0$ exponentially.  Thus, $\tilde{K}_{\eta i}(t)v$ decays to zero exponentially, and hence, $\lim_{t \rightarrow \infty} f_i (t) = 0$ exponentially.
Moreover, since $(A_i+ B_i K_{xi})$ is Schur,
by  Lemma \ref{lem0},  for any initial condition $\tilde{x}_i(0)$, $\lim_{t\rightarrow \infty} \tilde{x}_i(t) =0$ exponentially.
Thus, $\lim_{t\rightarrow \infty}\tilde{u}_i(t)=0$ by (\ref{sob2}),
and hence $\lim_{t\rightarrow \infty} {e}_i(t)=0$ by \eqref{sob0.2}.
\end{Proof}

Next, we will further consider solving Problem \ref{prob1}
by measurement output feedback control. Let $K_{xi}$ and
$K_{\eta i}(t)$ be defined as in the state feedback control law \eqref{ctrlu}.
Since $(C_{mi},A_i)$ is detectable, there exists $L_i$ such that
$A_i+L_iC_{mi}$ is Schur.
For $i=1,\dots,N$, we design the following measurement output feedback controller
\begin{subequations}\label{crl}
    \begin{align}
      u_i&=K_{xi}\xi_i+K_{\eta i}(t)\eta_i \label{crl.1}\\
      {\xi}_i (t+1) &= A_i\xi_i+B_iu_i+E_i\eta_i \nonumber \\
      &+L_i(C_{mi}\xi_i+D_{mi}u_i+F_{mi}\eta_i-y_{mi}). \label{crl.2}
      \end{align}
\end{subequations}

We have the following result.
\begin{Theorem}\label{thm2}
Under Assumptions \ref{ass1}-\ref{ass4},  let $\mu_1$,
 $\mu_2$, and $\mu_{3i}$ be the same as those in Theorem \ref{thm1}.  Then, Problem \ref{prob1} is solvable by the control law composed of
  \eqref{do} , \eqref{ctrlols} and \eqref{crl}.
\end{Theorem}
\begin{Proof}
Let $\hat{x}_i=\xi_i-x_i$. Then we have
  \begin{equation}
    \begin{aligned}
    &{\hat{x}}_i (t+1)=A_i\xi_i+B_iu_i+E_i\eta_i-A_ix_i-B_i u_i-E_iv+\\
      &L_i(C_{mi}\xi_i+D_{mi}u_i+F_{mi}\eta_i-C_{mi}x_i-D_{mi}u_i-F_{mi}v)\\
      &=A_i\hat{x}_i+E_i\tilde{\eta}_i+L_iC_{mi}\hat{x}_i+L_iF_{mi}\tilde{\eta}_i\\
      &=(A_i+L_iC_{mi})\hat{x}_i+(E_i+L_iF_{mi})\tilde{\eta}_i.
    \end{aligned}
  \end{equation}
  Since $\lim_{t\rightarrow\infty}\tilde{\eta}_i(t)=0$ exponentially and
  $(A_i+L_iC_{mi})$ is Schur, by Lemma \ref{lem0},
  $\lim_{t\rightarrow\infty}\hat{x}_i(t)=0$ exponentially.
  Note that in this case, \eqref{sob0.1} and \eqref{sob0.2} still hold.
  Next, similar to \eqref{sob2}, a simple calculation gives
  \begin{equation}
    \begin{aligned}
      \tilde{u}_i (t) &= K_{xi} (\tilde{x}_i+\hat{x}_i + X_{i}v) + K_{\eta i} (t) ( \tilde{\eta}_i  + v) - U_{i} v \\
      &=K_{xi}\tilde{x}_i +K_{xi}\hat{x}_i+ K_{\eta i}(t)\tilde{\eta}_i+\tilde{K}_{\eta i}(t)v.
    \end{aligned}
\end{equation}
  The rest of the proof follows from the proof of Theorem \ref{thm1} by noticing that $K_{xi}\hat{x}_i$ will also
  decay to zero exponentially.
  \end{Proof}

%\begin{Remark}
%  In comparison with the state feedback controller
%  given in \cite{shtac12}, which
%  is quoted as follows
%  \begin{equation}\label{ssfc}
%    u_i=K_{xi}x_i+K_{\eta i}\eta_i,
%  \end{equation}
%  we have
%  replaced the feedforward gain $K_{\eta i}$ by the
%  estimated feedforward gain $k_{\eta i}(t)$.
%  It is proved in Theorem \eqref{thm1} that the cooperative output regulation
%  problem can be solved by the state feedback controller \eqref{ctrlu}.
%\end{Remark}

\section{An Example}
\begin{figure}
\begin{center}
\scalebox{0.9}{\includegraphics[129,548][269,726]{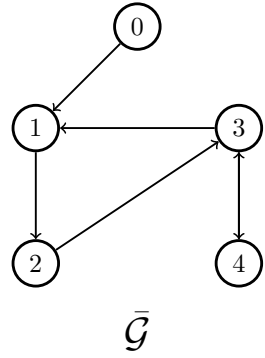}}
\caption{Communication Graph $\bar{\mathcal{G}}$.}\label{topo}
\end{center}
\end{figure}

In this section, we consider the cooperative output regulation problem of four agents.
The leader system is given by \eqref{2.2d} with
\begin{equation}
    S=\left[
         \begin{array}{cc}
           \cos \frac{\pi}{4} & \sin \frac{\pi}{4} \\
           - \sin \frac{\pi}{4} & \cos \frac{\pi}{4}
          \end{array}
       \right].
\end{equation}
Clearly,  Assumption \ref{ass1} is satisfied.
%With this $S$, the solution to \eqref{2.2d} with any initial condition $ v (0) = \mbox{col} (\alpha, \beta)$ is given by
%\EQQ
%v (t)= \left[
%                    \begin{array}{c}
%                      \alpha \\
%                      \beta \\
%                    \end{array}  \right] \cos \frac{\pi}{4} t +  \left[
%                    \begin{array}{c}
%                     \beta \\
%                      - \alpha \\
%                    \end{array}  \right] \sin \frac{\pi}{4} t
%                    \ENN

The four followers are given by \eqref{2.1} with
\begin{equation*}
    \begin{aligned}
     A_i&=\left[
                    \begin{array}{cc}
                      0 & 1 \\
                      0 & 0 \\
                    \end{array}
                  \right],B_i=\left[
                                \begin{array}{c}
                                  0 \\
                                  1 \\
                                \end{array}
                              \right],C_i=C_{mi}=\left[
                                \begin{array}{c}
                                  1 \\
                                  0 \\
                                \end{array}
                              \right]^T \\
   D_i & = D_{mi}=0, F_i = F_{mi}=[-1,0], E_i= \left[ \begin{array}{cc}
                                  0 & 2i-1\\
                                  0 & 1 \\
                                \end{array}\right]\\
   i &= 1, 2, 3, 4
    \end{aligned}
\end{equation*}
Thus,  Assumption \ref{ass2} is also satisfied.
Let us take $K_{xi} = [0.2, 0]$ such that, for $i = 1, 2, 3, 4$,  the eigenvalues of $(A_i+B_i K_{xi})$ are $\{0.447, -0.447\}$.

It can be verified that, for $i = 1, 2, 3, 4$,
the solution to the regulator equations are as follows:
\EQQ
X_i &= & \left[ \begin{array}{cc}
                                 1 & 0 \\
                                  \cos \frac{\pi}{4} & \sin \frac{\pi}{4} - (2i-1)
                                \end{array}\right],\\
U_i &=& [\cos \frac{\pi}{4},  \sin \frac{\pi}{4} - (2i-1) ] S - [0, 1] \\
i &= & 1, 2, 3, 4.
\ENN
Thus, Assumption \ref{ass3} is verified

The communication graph is shown in Fig. \ref{topo}. Thus, Assumption \ref{ass4} is satisfied with
\EQQ
H &=& \left[ \begin{array}{cccc}
2 & 0 &   -1 & 0 \\
                             -1 & 1 & 0 & 0 \\
                             0 & - 1 & 2  & -1 \\
                             0 & 0&  -1 & 1
                                \end{array}\right].
\ENN

By Theorems \ref{thm1} and \ref{thm2}, the cooperative output regulation problem for this system is solvable by both state and output feedback control laws.
%we only show the performance of state feedback control law here.
Since $\sigma (H) = \{2.420\pm j 0.606, 1, 0.161\}$,  we can take  $\mu_1$ such that $0 < \mu_1 < 0.778$.
Since $ \sigma (S) = 0.707 \pm j 0.707$, by (\ref{conmu2}),  for all $\mu_2$ satisfying $0 < \mu_2 < 1.414$, (\ref{obt}) is stable.
To obtain $\mu_{3i}$, note that, for $i = 1, 2, 3,4$,
$$
\sigma (Q_i^T Q_i) = \{0.199, 0.199, 1.56, 1.56, 3.25, 3.25\}
$$
Thus, we can take $\mu_{3i}$ such that $0 < \mu_{3i}
 < 0.615$.
Thus, we can design a state feedback control law with $\mu_1 = 0.3$, $\mu_2 = 0.4$, and, for $1, 2, 3, 4$, and $\mu_{3i} = 0.1$.
It can be verified that Assumption \ref{ass2.1} is also satisfied. Thus, it is also possible to obtain an output feedback control law.

\begin{figure}
\begin{center}
\scalebox{0.30}{\includegraphics[86,172][540,608]{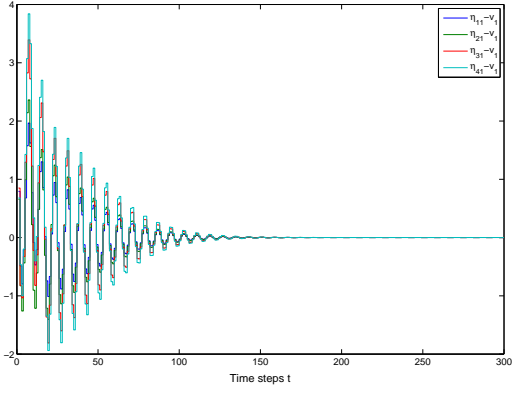}}
\vspace{0.3cm}
\caption{The estimation errors of the first component of the leader's signal by the adaptive distributed observer.}\label{a1}
\vspace*{0.3cm}
\end{center}
\end{figure}

\begin{figure}
\begin{center}
\scalebox{0.30}{\includegraphics[86,172][540,608]{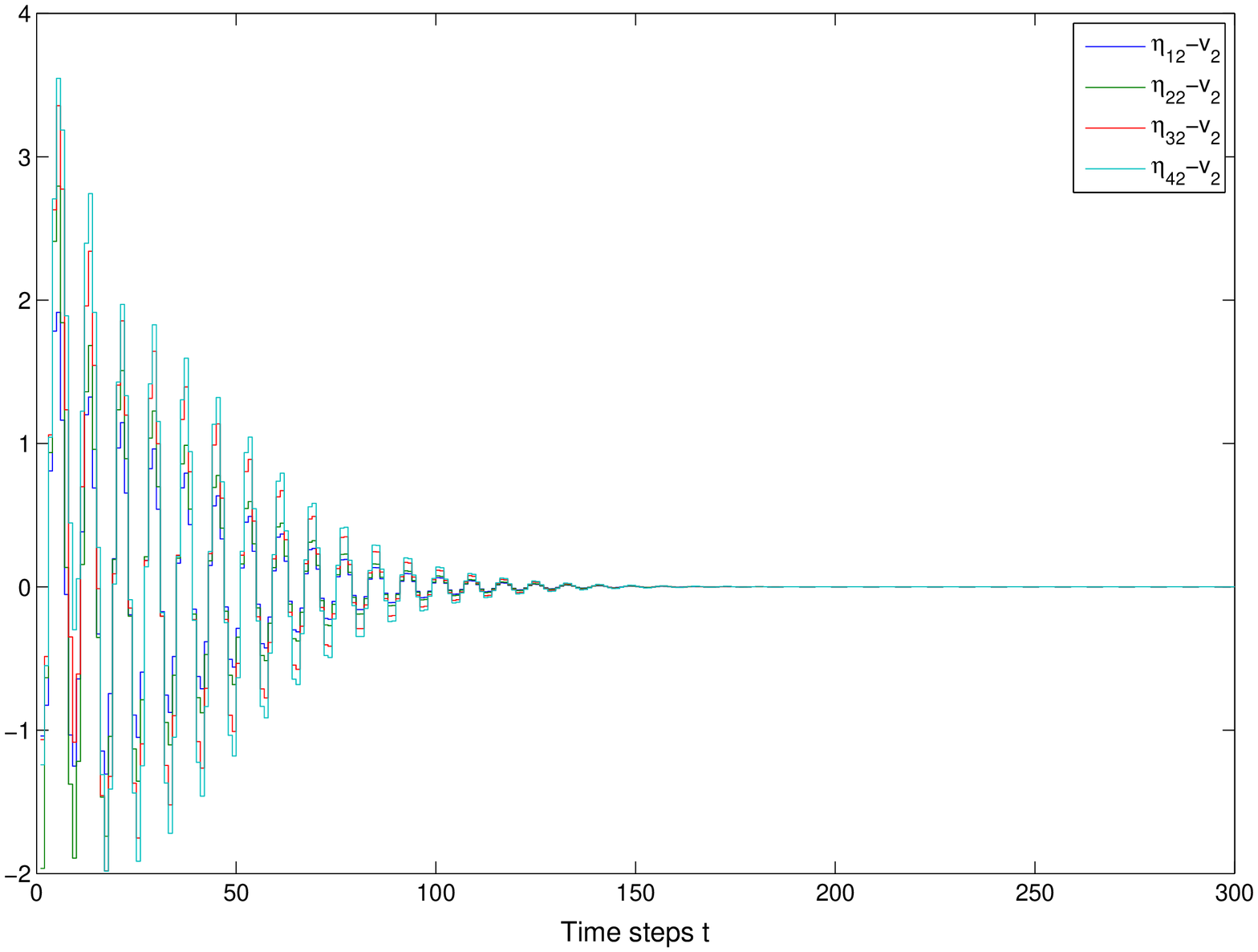}}
\vspace{0.3cm}
\caption{The estimation errors of the second component of the leader's signal by the adaptive distributed observer.}\label{a2}
\end{center}
\end{figure}

\begin{figure}
\begin{center}
\scalebox{0.30}{\includegraphics[86,172][540,608]{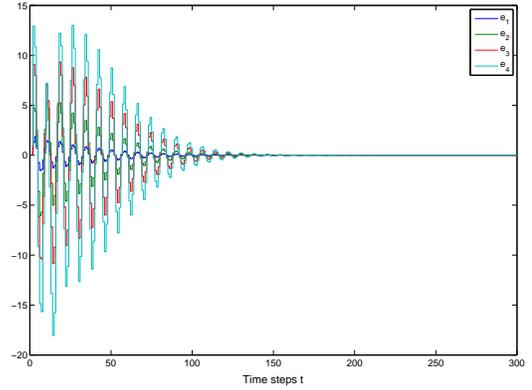}}
%\scalebox{0.5}{\includegraphics[94,184][527,596]{fig3.eps}}
\vspace{0.3cm}
\caption{The tracking errors of the four followers.}\label{a3}
\end{center}
\end{figure}

With the initial condition $v (0) = \mbox{col} (0, 2)$, the solution to the leader system is $v (t) = \mbox{col} (2 \sin \frac{\pi}{4} t, 2 \cos \frac{\pi}{4} t) $.
Figs. \ref{a1} and \ref{a2} show the estimation errors of the first and second components of the leader's signal by the adaptive distributed observer, respectively.
Fig. \ref{a3} further shows the tracking errors of the four followers under state feedback control law. Satisfactory tracking performance is observed.

%The observer gains are chosen as
%$\mu_1=\mu_2=10$, $\mu_3=20$. The initial values are
%given by $v(0)=[1,0.2]^T$, and for $i=1,2,3,4$, $x_i(0)=0$,
%$S_i(0)=0$, $\eta_i(0)=0$, $\zeta_i(0)=0$, $\xi_i(0)=0$. The regulated outputs

%\vspace*{0.3cm}

\section{Conclusion}
In this paper, we have studied the adaptive cooperative output regulation problem for discrete-time linear multi-agent
systems utilizing an adaptive distributed
observer.
%Both distributed state feedback and distributed measurement output feedback controllers are synthesized to
%solve the cooperative output regulation problem.
Compared with the existing distributed
observer based approach,  the approach of this paper does not require  that the
system matrix of the leader system be known by each follower.

One of the common challenges for the control of multi-agent systems is the the delay and the communication delay \cite{Hu1, olfatisaber2004, let}.
A natural extension of the current paper is to further consider same problem for systems with such delays.

\end{document}